# Modelling the Transmission Dynamics of Nipah Virus with Optimal Control

B. I. Omede[1], P. O. Ameh[1], A. Omame[2] [*] and B. Bolaji[1]

[1] Kogi State University, Anyigba; Nigeria.

[2] Federal University of Technology, Owerri; Nigeria

[*] Corresponding author E-mail address: omame2020@gmail.com, andrew.omame@futo.edu.ng

**Abstract**

A deterministic mathematical model is formulated and analyzed to study the transmission dynamics of Nipah virus both qualitatively and numerically. Existence and stability of equilibria were investigated and the model was rigorously analyzed. We then incorporated time dependent controls on the model, using Pontryagin's Maximum Principle to derive necessary conditions for the optimal control of the disease. We examined various combination strategies so as to investigate the impact of the controls on the spread of the disease. Through the numerical results, we found out that the optimal control strategies help to reduce the burden of the diseases significantly.

## 1.0   Introduction

Nipah Virus Infection (NiV), is a *zoonotic virus* (as it is transmitted from animals to humans) that causes outbreaks of fatal encephalitis in humans [1, 2, 4]. The virus derived its name from the fact that it was first isolated from a patient from Sungai Nipah village in Malaysia, thus it was christened Nipah according to the name of that village [1]. It was first recognized in a large outbreak of 276 reported cases in peninsular Malaysia and Singapore from September 1998 through May 1999 (see for examples [3-6, 11, 18]). Fruit bats that are migratory which are of the genus *Pteropus* have been identified as natural reservoirs of NiV [6]. Infected bats shed virus in their excretion and secretion such as saliva, urine, semen and excreta but they are symptomless carriers. The NiV is highly contagious among pigs, spread by coughing while direct contact with infected pigs was identified as the predominant mode of transmission in humans when it was first recognized in a large outbreak in Malaysia in 1999 [14]. It was discovered that ninety percent of the infected people in the 1998-1999 outbreaks were pig farmers or had contact with pigs. Drinking of fresh date palm sap, possibly contaminated by fruit bats (P. giganteus) during the winter season, may have been responsible for indirect transmission of Nipah virus to humans while there is strong evidence indicative of human- to-human transmission of NiV as discovered in Bangladesh in 2004 [18]. During the outbreak in Siliguri, 33 health workers and hospital visitors became ill after exposure to patients hospitalized with Nipah virus illness, suggesting nosocomial infection [17]. Symptoms of NiV infection in humans are similar to that of influenza such as fever and muscle pain, while in some cases, inflammation of the brain occurs leading to disorientation or coma [15]. Encephalitis may present as acute or late onset. The case fatality rate ranges from 9 to 75%, while the Incubation period of the disease ranges from 4 to 18 days [15]. Human-to-human transmission of NiV has been reported in recent outbreaks demonstrating a risk of transmission of the virus from infected patients to healthcare workers through contact with infected secretions, excretions, blood or tissues. However, there is no effective treatment for Nipah virus disease, but ribavarin may alleviate the symptoms of nausea, vomiting, and convulsions [20]. Treatment is mostly focused on by managing fever and the neurological symptoms. Severely ill individuals need to be hospitalized and may require the use of a ventilator. A vaccine is being developed. A recombinant sub-unit vaccine formulation protects against lethal Nipah virus challenge in cats [21]. Some researchers worked on the clinical study of Nipah Virus infection, including Goh [14] in his work discussed the clinical features of patients with Nipah virus encephalitis among the pig farmers in Malaysia; with major suggestion that the mental obtundation and neurologic signs indicate that Nipah virus probably has a predilection for the central nervous system. His conclusion is that the disease causes a rapidly progressive encephalitis with high mortality rate with features suggestive of involvement of brain stem. Others have made significant contributions to clinical studies of Nipah Virus infection, such as [12-14, 16-19, 22].

Researchers have undergone formulation of mathematical models such as the one in [23, 34-40] for the study of the transmission dynamics of different diseases and has helped reshaped understanding of these contagious diseases vis-à-vis procurement of measures to combat them adequately and effectively. There are few works on mathematical modelling of the transmission dynamics of Nipah Virus infection in literature, which include that of Biswas [29] who developed mathematical model comprising ordinary differential equations for the transmission dynamics of Nipah virus which he solved numerically and then analyse the behaviour of the disease dynamics with the use of optimal control strategy. [30] discussed the analysis of a mathematical model for the control and spread of NiV infection with vital dynamics (birth and death rates which are not equal). [29] and [30] in their works did not take into consideration the treatment class in their models, meanwhile infection can still take place in the treatment class [30]. This however, have been painstakingly taken into consideration in this study.

Consequently, in this study, the main focus is to formulate and rigorously analyze a new model for the human to human transmission dynamics of Nipah Virus; the model is rigorously analyzed qualitatively and numerically. We then reformulated the model by incorporating controls on key parameters in the model towards the effective control of the dynamics of its transmission with a view to reducing the burden of transmission of the disease.

## 2.0 Model Formulation

The total population of the people in the community at time $(t)$, denoted by $N(t)$ is split into the mutually exclusive compartments of Susceptible people $S(t)$, Exposed people $E(t)$, Infected people $I(t)$, Treatment class $T(t)$, and Recovered people $R(t)$, so that

$$N(t) = S(t) + E(t) + I(t) + T(t) + R(t)$$

The population of susceptible people $S(t)$ is generated by the birth of children at a rate $\pi$ and diminished by individuals who acquire infection following contact with infected individuals at the rate $\frac{\beta I}{N}$. This population is further reduced by natural death (at a rate $(\mu)$; natural death occur in all epidemiological compartments at this rate), so that;

$$\frac{dS}{dt} = \pi - \frac{\beta IS}{N} - \mu S$$

The population of exposed individuals $E(t)$ is generated by susceptible individuals who came in contact with those who are infected at the rate $\frac{\beta I}{N}$. This population is diminished by exposed individuals who progressed to infected class at the rate $\sigma$ and natural death, so that;

$$\frac{dE}{dt} = \frac{\beta IS}{N} - \sigma E - \mu E$$

The population of the infected class $I(t)$ is generated by exposed individuals who progressed to infected class at a rate $\sigma$ and the population is further increased by individuals who are infected in the treatment class at a rate $v$. this population is diminished by infected individuals who progressed to treatment class at the rate $\gamma$, population of individuals who recovered naturally as a result of strong antibodies at the rate $\varepsilon_1$, population of disease-induced death of individuals at the rate $\delta$ and natural death rate $\mu$, so that;

$$\frac{dI}{dt} = \sigma E + vT - (\gamma + \varepsilon_1 + \delta + \mu)I$$

The population of the treated class $T(t)$ is generated by infected individuals who progressed to treatment class at the rate $\gamma$. This population is decreased by individuals who recovered due to treatment at the rate $\alpha$ and the individuals infected in the treatment class at the rate $v$. This population is further diminished by natural recovery in the treated class at the rate $\varepsilon_2$, the disease-induced death in the treatment class at the rate $\theta\delta$ and natural death rate $\mu$, so that;

$$\frac{dT}{dt} = \gamma I - (\alpha + v + \varepsilon_2 + \theta\delta + \mu)T$$

The population of the recovered class $R(t)$ is generated by recovery due to treatment at the rate $\alpha$, natural recovery in the treated class and infected class at rate $\varepsilon_2$ $and$ $\varepsilon_1$ respectively. While the population of the class diminish by natural death rate $\mu$, so that we have:

$$\frac{dR}{dt} = \alpha T + \varepsilon_1 I + \varepsilon_2 T - \mu R$$

Based on the above formulations and assumptions, the following deterministic system of non-linear differential equations is the model for the transmission dynamics of Nipah Virus; the associated state variables, parameters and flow diagram of the model are depicted below:

$$\frac{dS}{dt} = \pi - \frac{\beta IS}{N} - \mu S$$

$$\frac{dE}{dt} = \frac{\beta IS}{N} - \sigma E - \mu E$$

$$\frac{dI}{dt} = \sigma E + vT - (\gamma + \varepsilon_1 + \delta + \mu)I \qquad (1)$$

$$\frac{dT}{dt} = \gamma I - (\alpha + v + \theta\delta + \varepsilon_2 + \mu)T$$

$$\frac{dR}{dt} = \alpha T + \varepsilon_1 I + \varepsilon_2 T - \mu R$$

We itemize some of the main assumptions made in the formulation of model (1) as follows:

1. The disease induced death rate in the infected and treatment class $(\delta)$ and natural death rate $(\mu)$ in all classes are the same.
2. Infection can take place in the treatment class [30].
3. Natural recovery from the infection can take place perhaps due to strong antibodies.
4. Natural recovery can take place in the treatment class.

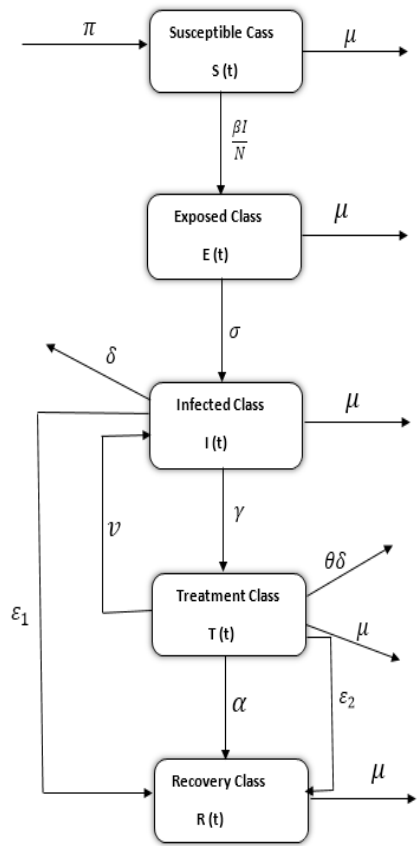

Figure 1. Flow diagram for the model

**Table 1** Description of variables and parameters of the Nipah Virus model (1)

| Variable | Description |
|---|---|
| $S$ | Population of susceptible individuals |
| $E$ | Population of exposed individuals |
| $I$ | Population of infected individuals |
| $T$ | Population of treated individuals |
| $R$ | Population of the recovered individuals |
| **Parameter** | **Description** |
| $\pi$ | Recruitment rate |
| $\beta$ | Transmission rate |
| $\sigma$ | Progression rate of people in exposed class to infectious class |
| $\gamma$ | Progression rate of people in infected class to treatment class |
| $\delta$ | Disease-induced death rate |
| $v$ | The rate of people infected in the treatment class |
| $\alpha$ | Recovery rate of infected people due to treatment |
| $\varepsilon_1$ and $\varepsilon_2$ | Natural recovery rate of people in infected and treatment class respectively |
| $\theta$ | Modification factor for progression to death by treated individuals |
| $\mu$ | Natural death rate |
| $N$ | Total number of the population |

## 2.1 Positivity of solution

Since the model (1) monitors human population, it is pertinent that all its state variables and associated parameters are non-negative for all time, *t*. Hence, the following non-negativity results holds for the state variables in model (1).

**Theorem 2.1**
Let the initial data be $\{(S(0), E(0), I(0), T(0), R(0)) \geq 0\} \epsilon R_+^5$. Then, the solution set $\{S(t), E(t), I(t), T(t), R(t)\}$ of the model equations with positive initial data, will remain positive for all time $t > 0$.

**Proof:**

From the first equation of the model

$$\frac{dS}{dt} = \pi - \frac{\beta IS}{N} - \mu S$$

$$\frac{dS}{dt} = \pi - \left(\frac{\beta I}{N} + \mu\right)S \geq -\left(\frac{\beta I}{N} + \mu\right)S$$

$$\frac{dS}{dt} \geq -\left(\frac{\beta I}{N} + \mu\right)S$$

By integrating both sides, we have:

$$\int \frac{dS}{S} \geq \int -\left(\frac{\beta I}{N} + \mu\right) dt$$

$$\ln S \geq -\left(\frac{\beta I}{N} + \mu\right) t + C_1$$

Taking the exponential of both sides, we have:

$$S(t) \geq e^{-\left(\frac{\beta I}{N}+\mu\right)t + C_1}$$

$$S(t) \geq C_1 e^{-\left(\frac{\beta I}{N}+\mu\right)t}$$

By applying the initial condition $t = 0$, $S(0) \geq C_1$, we have:

$$S(t) \geq S(0) e^{-\left(\frac{\beta I}{N}+\mu\right)t} \geq 0, \qquad \text{since } \left(\frac{\beta I}{N} + \mu\right) > 0$$

Similarly, it can be shown that other state variables: $S(t) > 0, E(t) > 0, I(t) > 0, T(t) > 0$, and $R(t) > 0$. Thus all solutions of the model (1) remain positive for all non-negative initial conditions.

## 2.2 Invariance Property

There is a need to prove that all the state variables of the model are non-negative for all time (t), for the model (1) to be epidemiologically meaningful. In other words, the solutions of the model (1) with positive initial data will remain positive for all $t \geq 0$; we do this below:

**Lemma 2.2** The region: $D = \{(S, E, I, T, R) \in R_+^5 : N \leq \frac{\pi}{\mu}\}$ is positively-invariant and attracts all the solution in $R_+^5$

**Proof**: By adding the equations in the model system (1), we have

$$\frac{dN}{dt} = \pi - (S + E + I + T + R)\mu - \delta I - \theta \delta T$$

i.e. $\frac{dN}{dt} = \pi - N\mu - (I + \theta T)\delta$

The above equation can be re-written as;

$$\frac{dN}{dt} \leq \pi - N\mu$$

$$\frac{dN}{dt} + \mu N \leq \pi$$

Finding the integrating factors, $(IF) = e^{\int \mu dt} = e^{\mu t}$ we have:

$$e^{\mu t}\frac{dN}{dt} + \mu N e^{\mu t} \leq \pi e^{\mu t}$$

$$\frac{d(Ne^{\mu t})}{dt} \leq \pi e^{\mu t}$$

Integrating on both sides we get:

$$Ne^{\mu t} \leq \frac{\pi}{\mu} e^{\mu t} + c$$

Where $c$ is a constant of integration, therefore,

$$N \leq \frac{\pi}{\mu} + ce^{-\mu t}$$

Using the initial conditions; where $t = 0, N(0) = N_0$

$$N_0 - \frac{\pi}{\mu} \leq c$$

$$N \leq \frac{\pi}{\mu} + \left(N_0 - \frac{\pi}{\mu}\right)e^{-\mu t}$$

By applying theorem given by [28] on the differential inequality above, we obtain $0 \leq N \leq \frac{\pi}{\mu}$ as $t \to \infty$. Thus, the region D is a positively-invariant set under the flow described by model (1) so that no solution path leaves through the boundary of region D. Hence, it is sufficient to consider the dynamics of the model (1) in region D. Thus, in this region, the model can be considered as been epidemiologically and mathematically well posed.

### 3.0 Stability of the model

A major requirement of an epidemiological model is that it is stable (local and global stability). Therefore, in this section we analyze the model (1) for its stability property.

### 3.1 TREATMENT-FREE MODEL

Before we analyze model (1), it is necessary to study a version of it in the absence of treatment which is obtained by setting $T = \theta = v = \gamma = \alpha = 0$ in model (1) given by:

$$\frac{dS}{dt} = \pi - \frac{\beta IS}{N} - \mu S$$

$$\frac{dE}{dt} = \frac{\beta IS}{N} - (\sigma + \mu)E$$

$$\frac{dI}{dt} = \sigma E - (\varepsilon_1 + \delta + \mu)I$$

$$\frac{dR}{dt} = \varepsilon_1 I - \mu R \tag{2}$$

For the above treatment-free model (2), it can be shown that the region:

$$D_1 = \left\{(S, E, I, R) \in R_+^4 : S + E + I + R \leq \frac{\pi}{\mu}\right\}$$

is positive-invariant and attracting, so that it is sufficient to consider the dynamics of the model (2) in $D_1$

### 3.1.1 Existence of Disease-Free Equilibrium (DFE)

The steady state where there is no infection (or there is absence of disease) a point where $E^* = I^* = R^* = 0$, is called disease-free equilibrium (DFE) which is given by:

$$\varepsilon_0 = (S^*, E^*, I^*, R^*) = \left(\frac{\pi}{\mu}, 0, 0, 0\right) \quad (3)$$

### 3.1.2 Local Stability of Disease-free equilibrium

The linear stability $\varepsilon_0$ can be established by using next generation operator method on system (2). By using the notation in [25], the matrices $F$ and $V$, for the new infection and the remaining transmission terms are respectively, given by:

$$F = \begin{pmatrix} 0 & \frac{\beta\pi}{N\mu} \\ 0 & 0 \end{pmatrix} \quad \text{and} \quad V = \begin{pmatrix} P_1 & 0 \\ -\sigma & P_2 \end{pmatrix}$$

Where $P_1 = \sigma + \mu$ and $P_2 = \varepsilon_1 + \delta + \mu$

It follows that reproduction number denoted by:

$$R_0 = \rho(FV^{-1})$$

Where $\rho$ is the spectral radius or largest eigenvalue of $FV^{-1}$.

Therefore, $R_0 = \dfrac{\beta\pi\sigma}{N\mu(P_1 P_2)}$

The *basic reproduction number* $R_0$ is a measurement of the potential for spreading disease in a population. Mathematically, $R_0$ is a threshold parameter for the stability of a disease-free equilibrium and is related to the peak and final size of an epidemic. It is defined as the expected number of secondary cases of infection which would occur due to a primary case in a completely susceptible population [23]. If $R_0 < 1$, then a few infected individuals introduced into a completely susceptible population will, on average, fail to replace themselves, and the disease will not spread. On the other hand, when $R_0 > 1$, then the number of infected individuals will increase with each generation and the disease will spread.

However, in many disease transmission models, the peak prevalence of infected hosts and the final size of the epidemic are increasing $R_0$, making it a useful measure of spread.

For the reproduction number as computed above to be meaningful, there is the need to interprete it, this will do as follows.

**Interpretation of the reproduction number for treatment-free model:**

There is a need to interpret our reproduction number:

$$R_0 = \frac{\beta \pi \sigma}{N\mu(P_1 P_2)}$$

To do this, we substitute into it the values of the variables $P_1, P_2$ therein as initially defined, to obtain:

$$R_0 = \frac{\beta \pi \sigma}{N\mu(\sigma+\mu)(\varepsilon_1+\delta+\mu)} \quad \text{that is} \quad R_0 = \frac{\beta \pi}{N\mu} \times \frac{\sigma}{\sigma+\mu} \times \frac{1}{\varepsilon_1+\delta+\mu}$$

Where $\frac{\beta \pi}{N\mu}$ is the transmission probability, $\frac{\sigma}{\sigma+\mu}$ is probability that an individual has become infectious (that is, moved from exposed to infected class) and $\frac{1}{\varepsilon_1+\delta+\mu}$ is the average time that a person is infectious.

To health care practitioners, the reproduction number is better simply and explicitly interpreted as follows:

$R_0 =$ [Infection rate at S] $\times$ [(fraction that has become infectious and move from exposed class E to infected class I)(duration in exposed class E)] $\times$ [duration in infected class I].

By using theorem 2 in [25], the following result is established.

**Theorem 3.1** The Disease-free Equilibrium (DFE) of the model (2), given by equation (3) above, is locally asymptotically stable (LAS) if $R_0 < 1,$ and unstable if $R_0 > 1$.

Biologically, the meaning of theorem 3.1 is that Nipah Virus can be eliminated from the community when the reproduction number is less than unity ($R_0 < 1$) if the initial sizes of the subpopulations of the model are in the basin of the attraction of $\varepsilon_0$.

It is necessary to show that the Disease-free Equilibrium (DFE) is globally asymptotically stable in the invariance region $D_1$. We establish this as follows.

### 3.1.3 Global Stability of Disease-free equilibrium

**Theorem 3.2** The DFE of the treatment-free model (2) is globally asymptotically stable if $R_0 \leq 1$ and unstable if $R_0 \geq 1$.

**Proof:** Consider the following Lyapunov function;

$$F = aE + bI$$

With Lyapunov derivative (where a dot represents differentiation with respect to time)

$$\dot{F} = a\dot{E} + b\dot{I}$$

With a little perturbation from the reproduction number we have; $a = \sigma$ and $b = P_1$.

Therefore

$$\dot{F} = \sigma\dot{E} + P_1\dot{I} = \sigma\left(\frac{\beta IS}{N} - P_1 E\right) + P_1(\sigma E - P_2 I)$$

$$= \frac{\beta IS\sigma}{N} - P_1\sigma E + P_1\sigma E - P_1 P_2 I = \frac{\beta IS\sigma}{N} - P_1 P_2 I = I\left(\frac{\beta S\sigma}{N} - P_1 P_2\right)$$

$$\dot{F} \leq P_1 P_2 I\left(\frac{\beta \pi \sigma}{N\mu P_1 P_2} - 1\right) \qquad \forall\, S = \frac{\pi}{\mu} \leq N$$

$$\dot{F} = P_1 P_2 I (R_0 - 1) \leq 0 \qquad \forall\, R_0 \leq 1$$

Since all the model parameters are non-negative, it follows that $\dot{F} \leq 0$ for $R_0 \leq 1$ with $\dot{F} = 0$ if and only if $E = I = 0$. Hence, $F$ is a Lyapunov function in the invariant region. Therefore, by LaSalle's invariance principle, every solution to the equations in the treatment-free model, with condition in the invariant region, approaches $\varepsilon_0$ as $t \to \infty$.

Thus, the above theorem shows that the classical epidemiological requirement of $R_0 \leq 1$ is the necessary and sufficient condition for the elimination of the disease from the community.

### 3.1.4 Endemic Equilibrium (EE)

The endemic equilibrium can be obtained when $(S, E, I, R) \neq 0$. To obtain the conditions necessary for the existence of endemic equilibrium point, the treatment-free model (2) is solved in terms of the force of infection, $\lambda = \frac{\beta I}{N}$

Solving model (2) at steady state gives:

$$S^{**} = \frac{\pi}{\lambda^{**} + \mu}, \quad E^{**} = \frac{\pi\lambda^{**}}{P_1(\lambda^{**} + \mu)}, \quad I^{**} = \frac{\sigma\pi\lambda^{**}}{P_1 P_2(\lambda^{**} + \mu)}, \quad R^{**} = \frac{\varepsilon_1 \sigma\pi\lambda^{**}}{\mu P_1 P_2(\lambda^{**} + \mu)},$$

Where $P_1 = \sigma + \mu$ and $P_2 = \varepsilon_1 + \delta + \mu$

From $\lambda^{**} = \dfrac{\beta I^{**}}{N}$  (4)

By Substituting the value of $I^{**}$ into (4) we have:

$$\lambda^{**} = \dfrac{\beta \sigma \pi \lambda^{**}}{N P_1 P_2 \left(\lambda^{**} + \mu\right)}$$

Dividing all through this by $\lambda^{**}$ and simplifying we obtain:

$$\lambda^{**} = \mu \left(\dfrac{\beta \sigma \pi - N \mu P_1 P_2}{N \mu P_1 P_2}\right) \quad \& \quad \lambda^{**} = \mu(R_0 - 1) \tag{5}$$

Clearly, from (5) for $\lambda^{**}$ to have a unique value, $R_0$ must be greater than one, that is $R_0 > 1$.

**Lemma 3.3** The Endemic equilibrium (EE) is local asymptotically stable if $R_0 > 1$ and unstable if $R_0 < 1$.

### 3.2    TREATMENT MODEL

#### 3.2.1   Existence of Disease-Free Equilibrium (DFE)

The steady state where there is no infection (or there is absence of disease) is called disease-free equilibrium. That is, a point where $E^* = 0, I^* = 0, T^* = 0$ and $R^* = 0$. The disease-free equilibrium (DFE) of the treatment model is given by:

$$\varepsilon_0^t = \left(S^*, E^*, I^*, T^*, R^*\right) = \left(\dfrac{\pi}{\mu}, 0, 0, 0, 0\right)$$

**Basic Reproduction Number**

Using the matrices $F$ and $V$, for the new infection and the remaining transmission terms respectively give:

$$F = \begin{pmatrix} 0 & \dfrac{\beta \pi}{N \mu} & 0 \\ 0 & 0 & 0 \\ 0 & 0 & 0 \end{pmatrix} \quad \text{and} \quad V = \begin{pmatrix} P_1 & 0 & 0 \\ -\sigma & P_2 & -v \\ 0 & -\gamma & P_3 \end{pmatrix}$$

Where $P_1 = \sigma + \mu, P_2 = \gamma + \varepsilon_1 + \delta + \mu$ and $P_3 = \alpha + v + \theta \delta + \varepsilon_2 + \mu$

It follows that reproduction number denoted by:

$$R_0^t = \rho(FV^{-1})$$

Where $\rho$ is the spectral radius or largest eigenvalue of $FV^{-1}$ is given by:

$$R_0^t = \frac{\beta\pi\sigma P_3}{N\mu P_1(P_2 P_3 - v\gamma)}$$

For the reproduction number as computed above to be meaningful, there is the need to interprete it, this we do as follows.

**Interpretation of the reproduction number for treatment model:**

To interpret the reproduction number we substitute into it the values of the variables $P_1, P_2$ and $P_3$ therein as initially defined, to obtain:

$$R_0^t = \frac{\beta\pi\sigma(\alpha + v + \theta\delta + \varepsilon_2 + \mu)}{N\mu(\sigma + \mu)((\gamma + \varepsilon_1 + \delta + \mu)(\alpha + v + \theta\delta + \varepsilon_2 + \mu) - v\gamma)}$$

$$R_0^t = \frac{\beta\pi}{N\mu} \times \frac{\sigma}{(\sigma + \mu)} \times \frac{(\alpha + v + \theta\delta + \varepsilon_2 + \mu)}{(\gamma + \varepsilon_1 + \delta + \mu)(\alpha + v + \theta\delta + \varepsilon_2 + \mu) - v\gamma}$$

$R_0^t =$ [Infection rate at S] × [(fraction that has become infectious and move from exposed to infected class)(duration in exposed class E)] × [(fraction that survives treatment class T and moved to Recovered class R)(duration in infected class I) (duration in treatment class T) difference those that are infected in the treatment class]

**Theorem 3.4** The DFE of the treatment model is locally asymptotically stable if $R_0^t < 1$ and unstable if $R_0^t > 1$

### 3.2.2 Global Stability

**Theorem 3.5** The DFE of the treatment model is globally asymptotically stable if $R_0^t \leq 1$ and unstable if $R_0^t \geq 1$.

**Proof:** To prove this we will consider the following Lyapunov function:

$$F = aE + bI + cT$$

With Lyapunov derivative (where a dot represents differentiation with respect to time)

$$\dot{F} = a\dot{E} + b\dot{I} + c\dot{T}$$

With a little perturbation from the reproduction number from the treatment model we have: $a = \sigma P_3, b = P_1 P_3$ and $c = vP_1$

$$\dot{F} = \sigma P_3 \dot{E} + P_1 P_3 \dot{I} + v P_1 \dot{T}$$

$$= \sigma P_3 \left( \frac{\beta IS}{N} - P_1 E \right) + P_1 P_3 (\sigma E + vT - P_2 I) + v P_1 (\gamma I - P_3 T)$$

$$= \frac{\beta IS \sigma P_3}{N} - P_1 P_3 \sigma E + P_1 P_3 \sigma E + P_1 P_3 vT - P_1 P_2 P_3 I + P_1 v \gamma I - P_1 P_3 vT$$

$$= \frac{\beta IS \sigma P_3}{N} - P_1 P_2 P_3 I + P_1 v \gamma I = \frac{\beta IS \sigma P_3}{N} - P_1 I (P_2 P_3 - v\gamma)$$

$$= \left( \frac{\beta S \sigma P_3}{N} - P_1 (P_2 P_3 - v\gamma) \right) I = P_1 (P_2 P_3 - v\gamma) I \left( \frac{\beta S \sigma P_3}{N P_1 (P_2 P_3 - v\gamma)} - 1 \right)$$

$$\leq P_1 (P_2 P_3 - v\gamma) I \left( \frac{\beta \pi \sigma P_3}{N \mu P_1 (P_2 P_3 - v\gamma)} - 1 \right) \qquad \forall\, S = \frac{\pi}{\mu} \leq N$$

$$\dot{F} = P_1 (P_2 P_3 - v\gamma) I (R_0^t - 1) \leq 0 \qquad \forall R_0^t \leq 1$$

Since all the model parameters are non-negative, it follows that $\dot{F} \leq 0$ for $R_0^t \leq 1$ with $\dot{F} = 0$ if and only if $E = I = T = 0$. Hence, $F$ is a Lyapunov function in the invariant region, therefore, by the LaSalle's invariance principle, every solution to the equations in the treatment model, with condition in the invariant region, approaches $\varepsilon_0^t$ as $t \to \infty$. Thus, the above theorem shows that the classical epidemiological requirement of $R_0^t \leq 1$ is the necessary and sufficient condition for the elimination of the disease from the community.

### 3.2.3 Endemic Equilibrium

The endemic equilibrium can be obtained when $(S, E, I, T, R) \neq 0$. To obtain the conditions necessary for the existence of endemic equilibrium point, the treatment model is solved in terms of the force infection. By solving model (1) at steady state gives:

$$S^{**} = \frac{\pi}{\lambda^{**} + \mu}, \qquad E^{**} = \frac{\pi \lambda^{**}}{P_1 (\lambda^{**} + \mu)}, \qquad I^{**} = \frac{P_3 \sigma \pi \lambda^{**}}{P_1 (\lambda^{**} + \mu)(P_2 P_3 - v\gamma)},$$

$$T^{**} = \frac{\pi \sigma \gamma \lambda^{**}}{P_1 (\lambda^{**} + \mu)(P_2 P_3 - v\gamma)}, \qquad R^{**} = \frac{\pi \sigma \lambda^{**} (\varepsilon_1 P_3 + \gamma P_4)}{\mu P_1 (\lambda^{**} + \mu)(P_2 P_3 - v\gamma)}$$

Where $P_1 = \sigma + \mu, P_2 = \gamma + \varepsilon_1 + \delta + \mu, P_3 = \alpha + v + \theta \delta + \varepsilon_2 + \mu$ and $P_4 = \alpha + \varepsilon_2$

By substituting the value of $I^{**}$ into (4) and simplifying, we have:

$$\lambda^{**} = \mu\left(\frac{\beta\pi\sigma P_3 - N\mu P_1(P_2 P_3 - v\gamma)}{N\mu P_1(P_2 P_3 - v\gamma)}\right) \quad \text{that is} \quad \lambda^{**} = \mu\left(R_0^t - 1\right) \tag{6}$$

Clearly, from (6), for $\lambda^{**}$ to have a unique value, $R_0^t$ must be greater than one, that is $R_0^t > 1$

**Lemma 3.6** The Endemic equilibrium (EE) is local asymptotically stable if $R_0^t > 1$ and unstable if $R_0^t < 1$.

## 4    Numerical Simulation of the model without control

The parameters in table 1 will be used in the Simulation. Some values assigned to the parameters were obtained from the literature. The objective of these simulations is to illustrate some of the theoretical results obtained in this study.

**Table 1.** Base line values of the parameters of the model with total population (N) estimated at 164,700,000.

| S/N | Parameters | Values | Source |
|---|---|---|---|
| 1 | $\pi$ | 6102 $year^{-1}$ | [30] |
| 2 | $\beta$ | 0.75 $year^{-1}$ | [29] |
| 3 | $\sigma$ | 0.60 $year^{-1}$ | [30] |
| 4 | $\gamma$ | 0.97 $year^{-1}$ | Inferred from [29] |
| 5 | $\delta$ | 0.76 $year^{-1}$ | [30] |
| 6 | $v$ | 0.89 $year^{-1}$ | Inferred from [29] |
| 7 | $\alpha$ | 0.09 $year^{-1}$ | [30] |
| 8 | $\varepsilon_1$ | 0.0054 $year^{-1}$ | Inferred from [30] |
| 9 | $\varepsilon_2$ | 0.0061 $year^{-1}$ | Inferred from [30] |
| 10 | $\theta$ | 0.51 $year^{-1}$ | Inferred from [29] |

| 11 | μ | 0.000038642 $year^{-1}$ | [30] |

**Remarks from the plots**

For the models without control, It is observed that the result gotten from the numerical simulation shows that reproduction number has a direct relationship with the rate of transmission (β) and the rate of infection in the treated class (ν), as would be expected, when β and ν increases, the reproduction number increases. Thus when β and ν decreases the reproduction number decreases. If rate of transmission (β) and the rate of infection in the treated class $(v)$ is reduced, Nipah virus can be eradicated from the community. Based on the results of this work, we recommend that more attention should be given to the individuals in the treatment class and standard precautions like the use of hand gloves, hand sanitizers should be placed on ground in hospitals, isolation centers, schools and homes when caring for patients infected with Nipah virus so as to control the transmission of this disease.

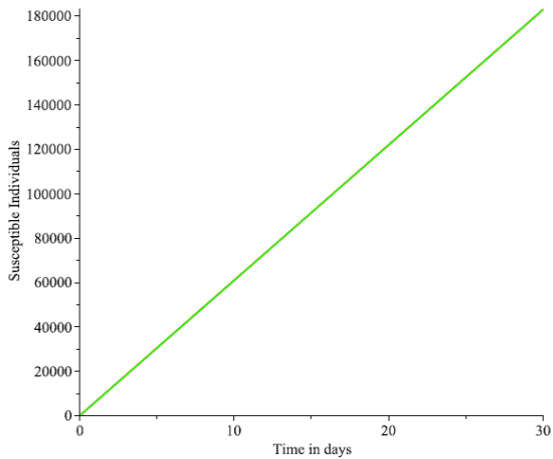
Figure 2. Simulation of Susceptible population

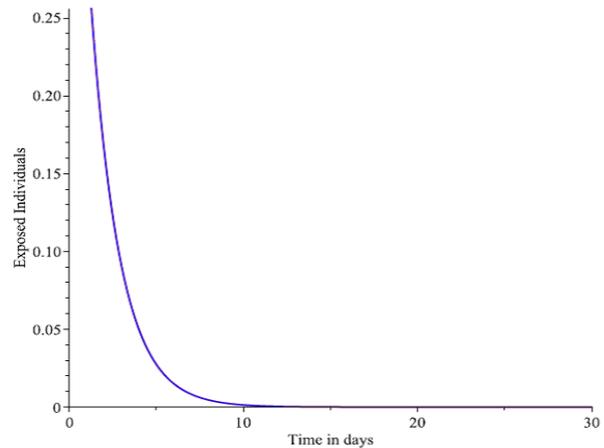
Figure 3. Simulation of Exposed individuals

The graph in figure 2 above shows the number of individuals susceptible to Nipah virus over the period of 30 days. As time increases the rate of susceptible individual increases as expected. In figure 3 above observe that as time increases, the number of exposed individuals decreases. That is, the curve of the Exposed class is flattened.

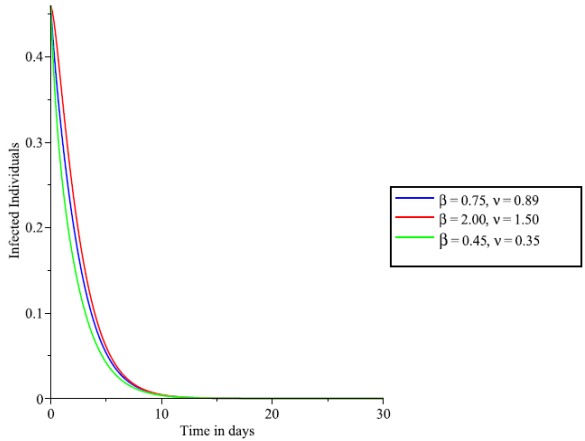 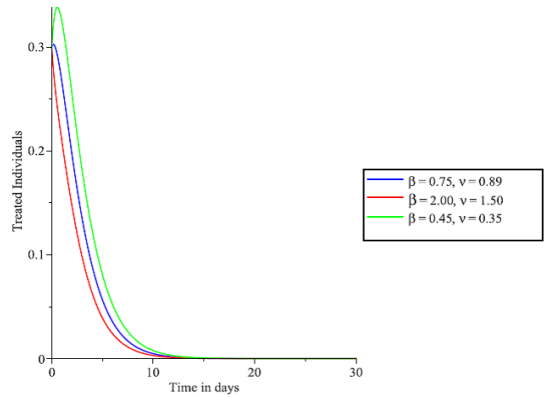

Figure 4. Simulation of infected

Figure 5. Simulation of Treated

In figure 4 above, observe that the infection level decreases with high recovery rate when the rate of infection in treated class and transmission rate reduces ($v = 0.35\ and\ \beta = 0.45$) as compared to when the rate of infection in the treated class and the rate of transmission is high ($v = 1.50\ and\ \beta = 2.00$).

From figure 5, observe that there is a decrease in the level of treated individuals when the rate of infection in the treated class and the rate of transmission is high ($v = 1.50\ \&\ \beta = 2.00$) due to the fact that individual died from the disease (via disease-induced death) before progressing to the treated class and the infected individual in the treated class will progress to infected class. As the rate of infection in the treated class and the rate of transmission reduces ($v = 0.35\ \&\ \beta = 0.45$), we witnessed an increased in the level of treated individuals.

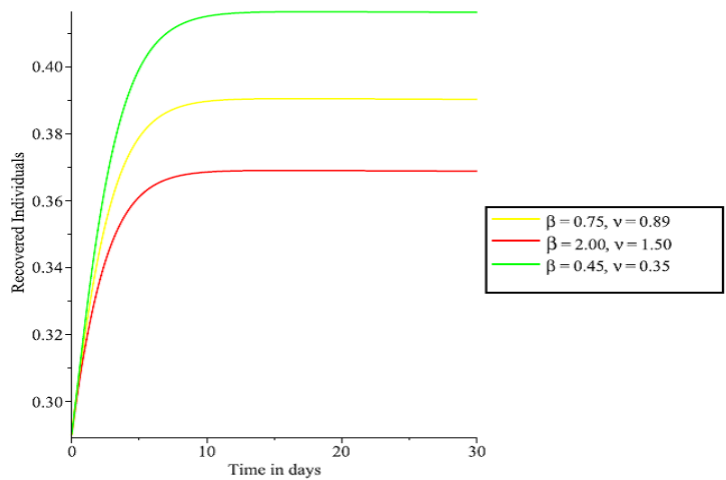

Figure 6. Simulation of Recovered individuals

From figure 6 above, observe that when the rate of infection in the treated class and the rate of transmission becomes high ($v = 1.50$ & $\beta = 2.00$), the chances of recovery from the disease reduces. While as the rate of infection in the treated class and the rate of transmission reduces ($v = 0.35$ & $\beta = 0.45$), the chances of recovery from Nipah virus increases.

## 5.0 Model with control strategy.

Optimal control theory is a concept that have been applied to various mathematical models for studying the transmission dynamics of contagious diseases, such as the ones discussed in [37-39]. It has no doubt assisted in no small way to gain insights into effective control of the spread of contagious diseases. Mentioning must be made of the fact that there is only very few works in literature on the mathematical modelling of Nipah Virus, apparently because, it is a re-emerging disease. [30] discussed the analysis of a mathematical model for the control and spread of Nipah Virus infection with vital dynamics, birth and death rates that are not equal and incorporated the quarantines of infectious individuals influenced by availability of isolation centers and surveillance coverage. They also develop a control model that can be used in studying how to reduce the spread and prevention of possible occurrence of Nipah Virus infections. They concluded that if there is increase in the values of the number of quarantined individuals as well as enhanced personal hygiene, then the spread of the disease will be reduced in a short possible time.

The model (1) is now modified with time dependent controls $u_1(t), u_2(t)$ and $u_3(t)$ where the dependent control strategies are as defined below:

$u_1(t)$: Prevention efforts to protect people in susceptible compartment from contacting Nipah virus.
$u_2(t)$: Treatment effort of Nipah virus for individuals in the Infectious compartment
$u_3(t)$: Treatment effort of Nipah virus for individuals in the Treatment compartment.

(7)

where the controls $u_1(t)$, $u_2(t)$ and $u_3(t)$ correspond to case finding and the coefficients $1-u_1(t), 1-u_2(t)$ and $1-u_3(t)$ represent the effort that prevents failure of the treatment in the drug-sensitive infectious individuals.

Using the controls above, the model with control strategy becomes:

$$\frac{dS}{dt} = \pi - (1-u_1)\frac{\beta IS}{N} - \mu S$$

$$\frac{dE}{dt} = (1-u_1)\frac{\beta IS}{N} - \sigma E - \mu E$$

$$\frac{dI}{dt} = \sigma E + \upsilon T - (\varepsilon_1 + \delta + \mu)I - (\gamma + u_2)I \tag{8}$$

$$\frac{dT}{dt} = (\gamma + u_2)I - (v + \vartheta\delta + \varepsilon_2 + \mu)T - (\alpha + u_3)T$$

$$\frac{dR}{dt} = (\alpha + u_3)T + \varepsilon_1 I + \varepsilon_2 T - \mu R$$

*The control parameters are Lebesque measurable set given by:*

$$\mu = \{u_1(t), u_2(t), u_3(t) : 0 \le u_1 \prec 1, 0 \le u_2 \prec 1, 0 \le u_3 \prec 1, 0 \le t \prec T\} \tag{9}$$

The objective function to be minimized is:

$$J = \min u_1, u_2, u_3 \int_0^t \left(b_1 E + b_2 I + \tfrac{1}{2} w_1 u_1^2 + \tfrac{1}{2} w_2 u_2^2 + \tfrac{1}{2} w_3 u_3^2\right) dt \tag{10}$$

In the objective function (10) above, the parameters $w_1, w_2$ and $w_3$ are the weight factors to help balance each term in the integrand in (10), so that none of the terms dominates. We explain the terms in the integrand in (10) above as follows:

(i) The term $b_1 E + b_2 I$ represents the cost associated with monitoring infected individuals in all stages.

(ii) The term $w_1 u_1^2$ and $w_2 u_2^2$ represents the cost associated with all forms of treatment for all infected individuals in the infectious classes.

(iii) The term $w_3 u_3^2$ represents the cost associated with all forms of treatment in the treatment class.

Our task here is to minimize the number of individuals in identified infectious classes and treatment classes giving them effective treatment in addition to reduction of their contact with the susceptible class such that we keep the cost of screening of immigrants low. Hence, we then minimize an objective function of a form that shows the trade-off needed in minimizing the number of infectious individuals vis-à-vis procuring of effective treatment for the infected individuals and the associated relevant cost of doing this. The associated relevant cost is made up of the cost of screening immigrants coming into the population, keeping them from contracting the disease and making them receive effective treatment.

Here we assume that the cost associated with the screening of immigrants coming into the population and the cost of keeping them from contracting the disease is nonlinear and take a quadratic form.

The goal is to find an optimal triplets $u_1(t)$, $u_2(t)$ and $u_3(t)$ such that:

$$J(u_1(t)^*, u_2(t)^*, u_3(t)^*) = \min J(u_1(t), u_2(t), u_3(t)) \text{ in } \Omega \tag{11}$$

Where $\Omega = \{(u_1(t), u_2(t), u_3(t)) \in L^1(0, t_f) \times L^1(0, t_f) | a_i \le u_i \le b_i\}$ and $a_i, b_i, i = 1, 2, 3$ are fixed non-negative constants.

*5.1 Analysis of optimal control model:* the necessary conditions that an optimal control triplet must satisfy is given by Pontryagin's Maximum Principle [32]. By using this principle we converts (8), (9) and (10) into a problem of minimizing a Hamiltonian $H$, point wisely with respect to the controls, $u_1(t)$, $u_2(t)$ and $u_3(t)$.

The Hamiltonian is given by:

$$H = b_1 E + b_2 I + \tfrac{1}{2} w_1 u_1^2 + \tfrac{1}{2} w_2 u_2^2 + \tfrac{1}{2} w_3 u_3^2$$
$$+ \lambda_1 \left[ \pi - (1-u_1)\frac{\beta IS}{N} - \mu S \right] + \lambda_2 \left[ (1-u_1)\frac{\beta IS}{N} - \sigma E - \mu E \right]$$
$$+ \lambda_3 [\sigma E + \upsilon T - (\varepsilon_1 + \delta + \mu)I - (\gamma + u_2)I] + \lambda_4 [(\gamma + u_2)I - (v + \vartheta\delta + \varepsilon_2 + \mu)T - (\alpha + u_3)T]$$
$$+ \lambda_5 [(\alpha + u_3)T + \varepsilon_1 I + \varepsilon_2 T - \mu R] \tag{12}$$

By applying the Pontryagin's Maximum Principle [32] and invoke the existence result for the optimal control triplets from [33], we obtain the following theorem:

**Theorem 5.1.** There exists an optimal control triplets $u_1(t)^*, u_2(t)^*$ and $u_3(t)$ and corresponding solution, $S^*, E^*, I^*, T^*$ and $R^*$, that minimizes $J(u_1(t), u_2(t), u_3(t))$ over $\Omega$. Furthermore, there exists adjoint functions, $\lambda_1, \lambda_2, \lambda_3, \lambda_4$ and $\lambda_5$ such that:

$$\lambda_1^1 = \frac{\beta I^* S^* (1-u_1)\lambda_1}{N^{*2}} + \frac{\beta I^* S^* (1-u_1)\lambda_2}{N^{*2}}, \qquad \lambda_2^1 = -1 - \frac{\beta I^* S^* (1-u_1)\lambda_1}{N^{*2}} + \frac{\beta I^* S^* (1-u_1)\lambda_2}{N^{*2}}$$

$$\lambda_3^1 = -1 - \frac{\beta I^* S^* (1-u_1)\lambda_1}{N^{*2}} + \frac{\beta I^* S^* (1-u_1)\lambda_2}{N^{*2}} - (\gamma + u_2)\lambda_4 + (\gamma + \delta + \mu + u_2 + \varepsilon_1)\lambda_5$$

$$\lambda_4^1 = -\frac{\beta I^* S^* (1-u_1)\lambda_1}{N^{*2}} + \frac{\beta I^* S^* (1-u_1)\lambda_2}{N^{*2}} - (\alpha + \theta\delta + \mu + v + \mu - \varepsilon_2)\lambda_4 + (\alpha + u_3 + \varepsilon_2)\lambda_5$$

$$\lambda_5^1 = \frac{\beta I^* S^* (1-u_1)\lambda_1}{N^{*2}} + \frac{\beta I^* S^* (1-u_1)\lambda_2}{N^{*2}} \tag{13}$$

with transversality conditions: $\lambda_i(t_f) = 0$, $i = 1, 2, ..., 5$ \hfill (14)

and $N^* = S^* + E^* + I^* + T^* + R^*$

The following characterization holds:

$$u_1(t)^* = \min\left[ \max\left( a_1, \frac{1}{B_1}(\lambda_2 p + \lambda_3 q - \lambda_1(p+q)), b_1 \right) \right]$$

$$u_2(t)^* = \min\left[ \max\left( a_2, \frac{1}{B_2}(\lambda_3 - \lambda_4)nI^*, b_2 \right) \right]$$

$$u_3(t)^* = \min\left[\max\left(a_3, \frac{1}{B_3}(\lambda_4 - \lambda_5)mT^*, b_3\right)\right] \quad (15)$$

**Proposition 3.1.** By virtue of Corollary 4.1 in [32], the convexity of the integrand of G with respect to $(u_1(t), u_2(t), u_3(t))$ guarantees the existence of an optimal pair, a priori boundedness of the state solutions, and the Lipschitz property of the state system with respect to the state variables. The proof of this proposition is left as *future work*.

**Proof of Theorem 5.1** We apply the Pontryagin's max principle to have:

$$\frac{d\lambda_1}{dt} = -\frac{\partial H}{\partial S}, \quad \lambda_1(t_f) = 0 \qquad \frac{d\lambda_2}{dt} = -\frac{\partial H}{\partial E}, \quad \lambda_2(t_f) = 0 \qquad \frac{d\lambda_3}{dt} = -\frac{\partial H}{\partial I}, \quad \lambda_3(t_f) = 0$$

$$\frac{d\lambda_4}{dt} = -\frac{\partial H}{\partial T}, \quad \lambda_4(t_f) = 0 \qquad \frac{d\lambda_5}{dt} = -\frac{\partial H}{\partial R}, \quad \lambda_5(t_f) = 0 \quad (16)$$

We shall also consider the optimality conditions:

$$\frac{\partial H}{\partial u_1} = 0, \qquad \frac{\partial H}{\partial u_2} = 0, \qquad \frac{\partial H}{\partial u_3} = 0,$$

and solving for $u_1(t)^*, u_2(t)^*$ and $u_3(t)^*$ subject to the state variables, the characterizations in (15) can be obtained, by taking into account the bounds on the control. By doing this, we have that:

$$\frac{\partial H}{\partial u_2} = 0, \qquad \frac{\partial H}{\partial u_3} = 0,$$

$$\frac{\partial H}{\partial u_1} = u_1 w_1 + \frac{\beta I^* S^* \lambda_1}{N^{*2}} - \frac{\beta I^* S^* \lambda_2}{N^{*2}} = 0, \qquad \Rightarrow u_1 = \frac{1}{w_1}\left(\frac{\beta I^* S^* \lambda_2}{N^{*2}} - \frac{\beta I^* S^* \lambda_1}{N^{*2}}\right)$$

$$\frac{\partial H}{\partial u_1} = u_2 w_2 - I^* \lambda_3 + I^* \lambda_4 = 0, \qquad \Rightarrow u_2 = \frac{1}{w_1}(I^* \lambda_3 - I^* \lambda_4)$$

$$\frac{\partial H}{\partial u_3} = u_3 w_3 - T^* \lambda_4 + T^* \lambda_5 = 0, \qquad \Rightarrow u_3 = \frac{1}{w_3}(I^* \lambda_4 - I^* \lambda_5) \quad (13)$$

It is observed that the optimality conditions obtained by taking the derivatives of the Hamiltonian (12) with respect to the controls only hold in the interior of the control set.
This ends the proof.

*5.2 Numerical Simulations of Optimal Control Model*

We obtained the optimality screening and treatment strategy by solving the optimality system which is made up of ten ODEs from the state and adjoint equations (a combination of 5 equations from the model with controls (8) and five equations from the adjoint equations (13)). We solved the optimality system by using an iterative method with a Runge-Kutta fourth order scheme. We solved the state system with an initial condition forward in time using MATLAB ODE 45, with a guess for the controls over the simulated time, while we solved the adjoint system with values at final time $t_f$ backward in time by using the current iteration solution of the state equations by using the same MATLAB ODE 45 routine. We then update the controls by using a convex combination of the previous controls and the value from the characterizations (15). This process is repeated and we stopped the process and iterations if the values of the unknowns at the previous iteration are very close to the ones at the present iteration (that is, until the control variables converges).

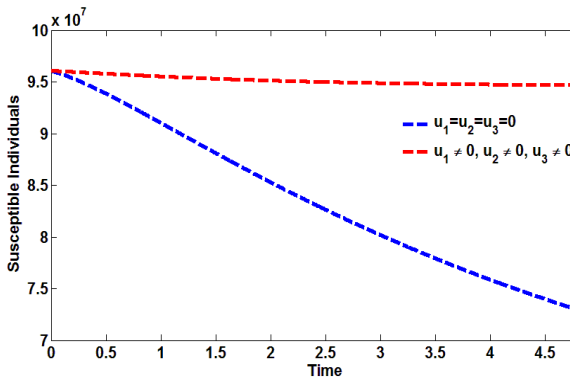

**Figure 5.1** Plot of Susceptible

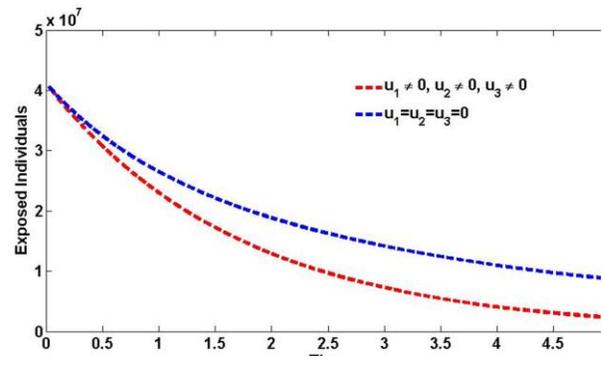

**Figure 5.2** Plot of Exposed

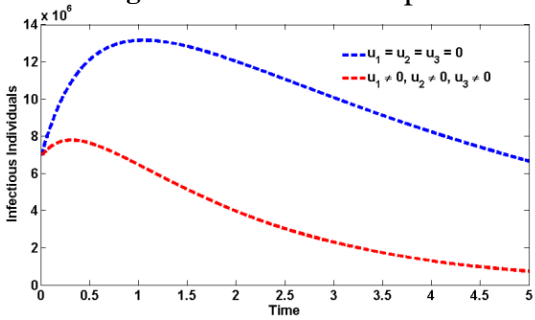

**Figure 5.3** Plot of Infectious

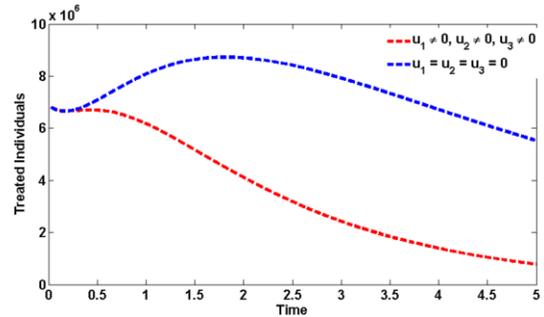

**Figure 5.4** Plot of Treated

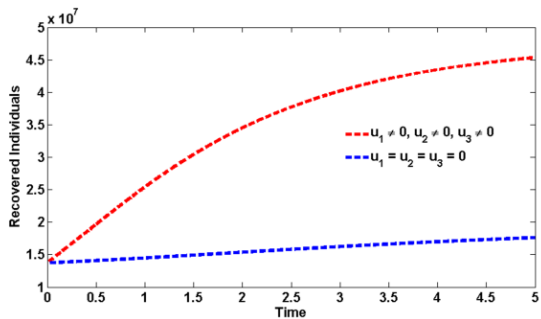

**Figure 5.5** Plot of Recovered Individuals

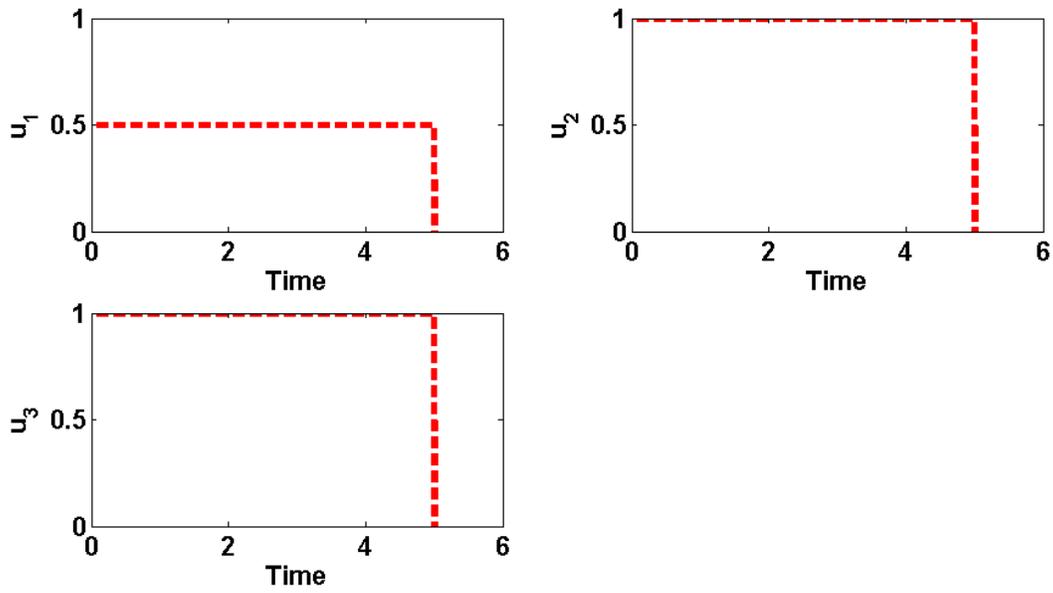

**Figure 5.6** Plot of Controls

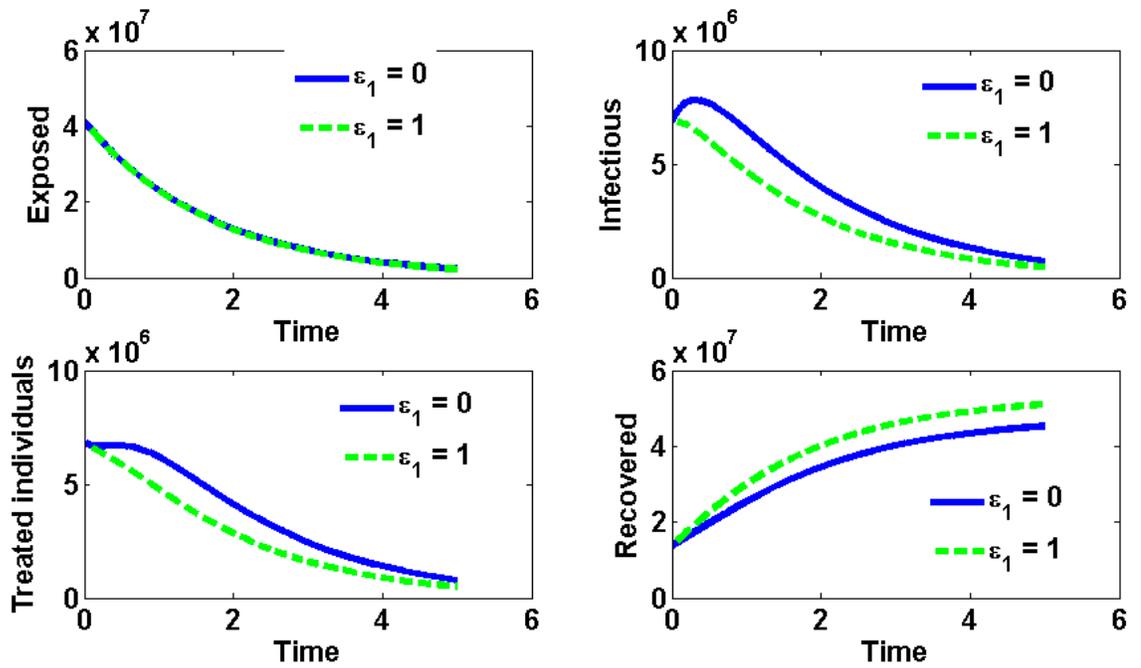

**Figure 5.7** Varying Epsilon 1 Plots with Controls

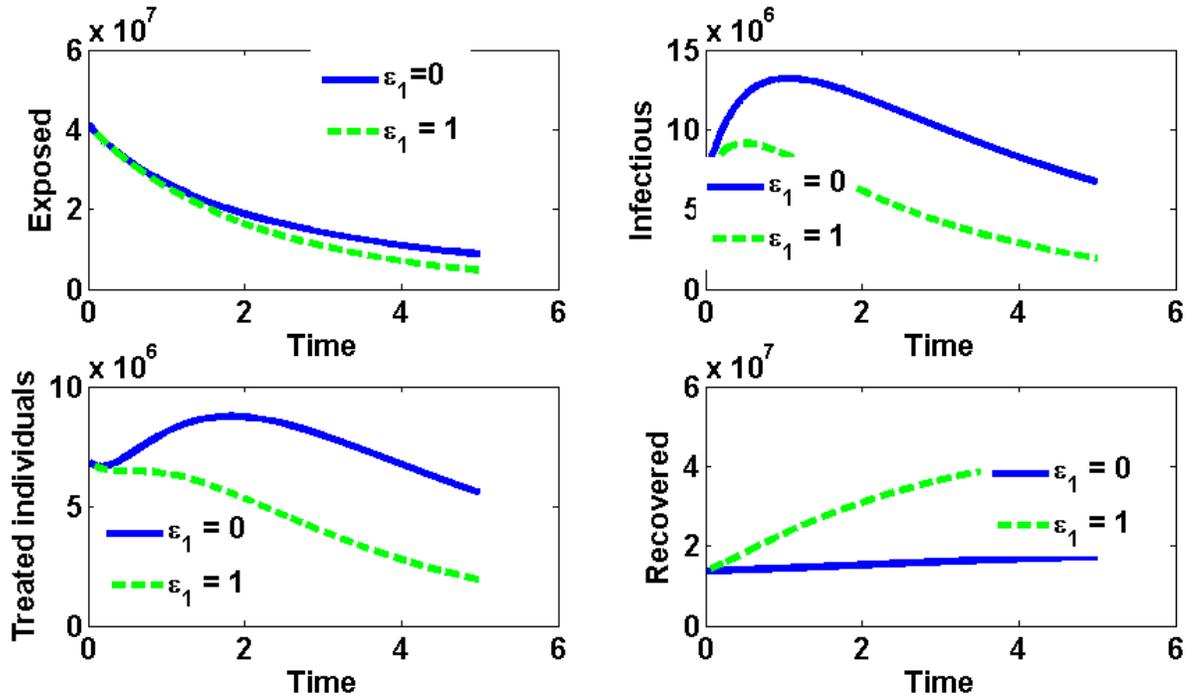

**Figure 5.8** Varying Epsilon 1 Plots without Controls

We plotted the epidemiological classes of the model with controls as a function of time with the controls; these is as shown in Figures 5.1, 5.2, …5.5. In Figure 5.1, we simulate the total number of Susceptible humans over time, in the presence of optimal control strategies. It is observed that the total number of susceptible individuals is more when controls are applied than when controls are not applied as would be expected. It is also seen in Figure 5.2, that the total number of individuals exposed to Nipah virus is less when controls are applied than when controls are not applied. The simulations of the total number of infectious and treated individuals, depicted in Figures 5.3 and 5.4, respectively show that the controls have significant impact in reducing the burden of the disease in a population. Figure 5.5 shows the simulation of the total number of recovered individuals over time, in the presence or absence of controls. It is observed that the total number of Recovered individuals is more when controls are applied than when controls are not implemented. The solution profiles of the three optimal controls $(u_1(t), u_2(t), u_3(t))$ are depicted in figure 5.6.

In Figures 5.7 and 5.8, it is shown that the various populations at different treatment rates for Infectious and treated individuals, when controls are applied lead to significant reduction in the burden of the disease. It is observed that the treatment rate for treated individuals has more impact on the dynamics of the model than the treatment rate for infectious, when controls are applied.

## 6.0   Discussion and Conclusion

In this work, we formulated a model for the transmission dynamics of Nipah Virus; both qualitative and numerical analysis of the model with (and without) treatment were done, the qualitative analysis shows that the models have two equilibria, the disease free equilibrium (DFE) which is locally asymptotically stable when a certain epidemiological threshold ($R_0$) is less than unity and that the disease will persist in the population when the threshold exceeds unity. Numerical results of the model without control shows that, the rate of transmission (β) and the rate of infection in the treated class (ν) has a direct relationship and great effect on the reproduction number.

In literature, there were very few works on the transmission dynamics of Nipah virus with control, consequently, we were motivated to contribute to knowledge in this area. We reformulated the model by incorporating controls into it. Three controls were incorporated viz: prevention efforts to protect people in susceptible compartment from contacting Nipah virus; treatment effort of Nipah virus for individuals in the Infectious compartment and treatment effort of Nipah virus for individuals in the Treatment compartment. From the numerical results of the model with control, we obtained the following:

(i) Treatment has a positive impact in reducing the burden of disease of Nipah virus.
(ii) The total number of infected individuals in the population is reduced by the application of time dependent controls.
(iii) The strategy which involves the implementation of all the controls is a good control strategy for Nipah Virus disease.

**References**

[1] Anno: Outbreak of Hendra-like virus – Malaysia and Singapore, 1998-99. Morb. Mort.

   Weekly Rep., 48 (13), 265-269. (1999).

[2] International Committee on Taxonomy of Viruses [database on the Internet]. Available from

   http://www.ncbi.nlm.nih.gov/ICTVdb/Ictv/index.htm. [Assessed August 2018].

[3] Yob, J. M, Field, H., Rashdi, A. M., Morrissy, C., van der Heide, B., Rota, P., et al: Nipah virus

   Infection in bats (order *Chiroptera*) in peninsular Malaysia. Emerg Infect Dis. 7:439–41.

   (2001).

[4] Chua, K. B, Bellini, W. J, Rota, P. A., Harcourt, B. H, et al.: Nipah virus: a recently emergent

   deadly paramyxovirus. Science. 288:1432-1435. (2000).

[5] Chua, K. B, Koh, C. L, Hooi, P. S, *et al:* Isolation of Nipah virus from Malaysian Island

   flying-foxes. *Microbes Infect.* 4 (2): 145–51. (2002).

   http://linkinghub.elsevier.com/retrieve/pii/S1286457901015222.

[6] Eaton, B. T, Broder, C. C, Middleton, D., Wang, L. F: Hendra and Nipah viruses: different and

dangerous. Nat. Rev. Microbiol. 4. 23–35. (2006).

[7] Hsu, V. P, Hossain, M. J, Parashar, U. D et al: Nipah virus encephalitis reemergence in

Bangladesh. Emerg. Infect. Dis. 10. 2082-2087. (2004).

[8] Reynes, J. M, Counor, D., Ong, S, Faure, C, Seng, V., Molia, S et al: Nipah virus in Lyle's

flyingfoxes, Cambodia. Emerg Infect Dis. 7. (2005). Available from

http://www.cdc.gov/ncidod/EID/vol11no07/04-1350.htm.

[9] Wacharapluesadee, S., Lumlertdacha., B., Boongird, K., Wanghongsa, S., Chanhome, L.,

Rollin, P. et al: Bat Nipah virus, Thailand. Emerg Infect Dis. 11:1949-51. (2005).

[10] Heymann, D. L: Henipavirus: Hendra and Nipah viral diseases. Control of communicable

Diseases Manua. 19th Edition. American Public Health Association. 275-278. (2008).

[11] Lehlé, C., Razafitrimo, G., Razainirina, J., *et al:* Henipavirus and Tioman virus antibodies

in pteropodid bats, Madagascar. *Emerging Infect. Dis*. 13 (1): 159– 61. *(*2007).

PMID 17370536. http://www.cdc.gov/ncidod/EID/13/1/159.htm.

[12] Hayman, D, Suu-Ire, R., Breed, A *et al:* "Evidence of henipavirus infection in West African

fruit bats". *PLOS ONE* 3 (7): 2739-2747. (2008).

http://www.plosone.org/article/info:doi/10.1371/journal.pone.0002739

[13] Mohd, M. N., Gan, C. H., Ong, B. L: Nipah virus infection of pigs in peninsular Malaysia

Rev. Sci. Tech. Off. Int. Epiz. 19:160-65. (2000).

[14] Goh, K. J, Tan, C. T, Chew, N.K, Tan, P. S. K., Kamarulzaman, A., Sarji, S. A., Wong, K.

T., Abdullah, B. J. J., Chua, K. B., Lam, S. K: Clinical features of Nipah virus Encephalitis

among pig farmers in Malaysia. N. Engl. J. Med. 342, 1229–1235. (2000).

[15] Hendra virus, ecological change and a new threat. http://scienceinpublic.com.au/

factsheets.htm. [Assessed August 2018]

[16] Luby, S. P., Rahman, M., Hossain, M.J., Blum, L.S., et al: Foodborne transmission of Nipah

virus, Bangladesh. Emerg. Infect. Dis. 12:1888-1894. (2006).

[17] Chadha, M. S, Comer, J. A., Lowe, L., Rota, P. A., Rollin, P. E., Bellini, W. J., et al: Nipah virus-associated encephalitis outbreak, Siliguri, India. Emerg. Infect. Dis. 12:235-40. (2006).

[18] Gurley, E., Montgomery, J. M., Hossain, M. J., Bell, M., Azad, A. K, Islam, M. R., et al: Person-toperson transmission of Nipah virus in a Bangladeshi community. Emerg. Infect. Dis. 13:1031-1037. (2007).

[19] Daniels, P., Ksiazek, T. G., Eaton, B. T.: Laboratory diagnosis of Nipah virus and Hendra infections. Microbes Infect. 3:289-295. (2001).

[20] Chong, H. T, Kamarulzaman, A., Tan, C.T, Goh, K. J, Thayaparan, T., Kunjapan, S. R, et al: Treatment of acute Nipah encephalitis with ribavirin. Ann Neurol. 49:810-813. (2001).

[21] McEachern, J. A., Bingham J., Crameri, G., Green, D. J., Hancock, T. J., Middleton D., et al: A recombinant subunit vaccine formulation protects against lethal Nipah virus challenge in cats. Vaccine. Volume 26, Issue 31, 3842-3852. (2008).

[22] Weingartl et al : Recombinant Nipah virus vaccines protect pigs against challenge. Journal of Virology. 80. 7929-38. (2006).

[23] Diekmann, O. and Heesterbeek, J. A. P.: Mathematical Epidemiology of Infectious Diseases. Model Building, Analysis and Interpretation. Wiley Series in Mathematical and Computational Biology. John Wiley and Sons, La.,Chichester. (2000).

[24] Murray, J. D : Mathematical Biology, an Introduction, Volume 17, Springer, New York. (2002).

[25] Chong, H. T., Jahangir Hossain, M., & Tan, C. T.: Differences in epidemiologic and clinical features of Nipah virus encephalitis between the Malaysian and Bangladesh outbreaks. *Neurology Asia*, *13*, 23 – 26. (2008).

[26] David, T. S. and Nicholas Johnson: The Role of Animals in Emerging Viral Diseases. (2014).

http://dx.doi.org/10.1016/B978-0-12-405191-1.00011-9

[27] Mills, J. N., Alim, A. N. M., Bunning, M. L., Lee O. B., Wagoner, K. D., Amman, B. R., et Al. Nipah virus infection in Malaysia. Emerg Infect Dis. (2009). Available from http://www.cdc.gov/EID/content/15/6/950.htm

[28] Birkhoff, G. and Rota, G. C: Ordinary Differential Equations. Ginn. (1982).

[29] Haider Ali Biswas, Mohidul Haque and Gomanth Duvvuru: A Mathematical Model for understanding the Spread of Nipah Fever Epidemic in Bangladesh. Emerg. Infect. Dis. (2015).

[30] Mitun Kumar Mondal, Muhammad Hanif & Md. Haider Ali Biswas: A mathematical analysis for controlling the spread of Nipah virus infection. International Journal of Modelling and Simulation. (2017). Available on http://dx.doi.org/10.1080/02286203.2017.1320820

[31] Van den Driessche, P. and Watmough, J.: Reproduction numbers and sub-threshold endemic equilibria for compartmental models of disease transmission. Math. Biosci. 180, 29 – 48. (2002).

[32] Pontryagrin, L.S., Boltyanskii, V.G., Gamkrelidze, R.V., et al.: The Mathematical theory of Optimal Processes. Wiley, New York. (1962).

[33] Fleming, W. H., Rishel R.W: Deterministic and Stochastic Optimal Control. Springer-Verlag, New York. (1975).

[34] Omame, A., Umana, R. A., Okuonghae, D., Inyama, S. C., Mathematical analysis of a two-sex Human Papillomavirus (HPV) model, International Journal of Biomathematics, 11(7)(2018).

[35] Omame, A., Okuonghae, D., Umana, R. A., Inyama, S. C., Analysis of a co-infection model for HPV-TB, Applied Mathematical Modelling, 77, 881-901 (2020).

[36] Omame, A., Okuonghae, D., Inyama, S. C., A mathematical study of a model for HPV with two high risk strains, in \textit{Mathematics Applied to Engineering, Modelling, and Social Issues} Studies in Systems, Decision and Control; F. Smith, H. Dutta and J. N. Mordeson (eds.) 200 (2020).

[37] Zaman, G., Kang, Y.H., Cho, G., Jung, I. H., Optimal strategy of vaccination & treatment in an SIR epidemic model, Math. Comput. Simulation (2016).

[38] Agusto, F.B., ELmojtaba, I.M., Optimal control and cost-effective analysis of malaria/ visceral leishmaniasis co-infection. PLoS ONE, 12 (2) (2017).

[39] Uwakwe, J.I., Inyama, S.C., Omame, A., Mathematical Model and Optimal Control of New-

Castle Disease (ND). Applied and Computational Mathematics. (2020) 9 (3) 70-84. doi: 10.11648/j.acm.20200903.14

[40] Umana, R.A., Omame, A., Inyama, S.C., Deterministic and Stochastic Models of the Dynamics of Drug Resistant Tuberculosis; FUTO Journals Series, (2016) 2 (2) 173-194.